 \documentclass[10pt,conference,a4paper]{ieeeconf} 
\usepackage{subfig}
\usepackage{float}
\usepackage{transparent}
\usepackage{amsmath}
\usepackage{amssymb}
\usepackage{graphicx}
\usepackage{color}
\usepackage{cite}
\usepackage{balance}

\usepackage{mathtools}

\mathtoolsset{showonlyrefs}

\def\<{\leqslant}           
\def\>{\geqslant}           

\def\[{[\![}
\def\]{]\!]}

\usepackage{transparent}
\usepackage{amsmath}
\usepackage{amssymb}
\usepackage{graphicx}
\usepackage{color}
\usepackage{fullpage}

\newcommand{\cZ}{{\cal Z}}

\newcommand{\gamovertwo} {{\frac{\gamma^2}{2}}}
\newcommand{\gammovertwo} \gamovertwo
\newcommand{\gammuovertwo} \gamovertwo
\newcommand{\gammuepstovertwo} \gamovertwo
\newcommand{\gammubarepstovertwo} \gamovertwo

\newcommand{\trace}{{\mbox{tr}}}

\newcommand{\rhatini}{{\hat{R}_0}}

\newcommand{\noncr} {\nonumber\\}

\newcommand{\beasnum}{\begin{eqnarray}}
\newcommand{\eeasnum}{\end{eqnarray}}
\newcommand{\beas}{\begin{eqnarray*}}
\newcommand{\eeas}{\end{eqnarray*}}

\newcommand{\be}{\begin{equation}}
\newcommand{\ee}{\end{equation}}
\newcommand{\ba}{\begin{array}}
\newcommand{\ea}{\end{array}}

\newtheorem{theorem}            {Theorem}[section]

\newtheorem{sideremark}         [theorem]{Remark}
\newtheorem{sideeg}           [theorem]{Example}
\newtheorem{sideconj}           [theorem]{Conjecture}

\def\argmin                     {\mathop{\rm argmin}}

\def \aeq {&=}
\def \adefeq {&:=}
\newcommand{\qed} {\hskip 0.2em\lower 0.7ex\hbox{\vbox{\hrule
\hbox{\vrule height 1.2ex\hskip 0.4em\vrule height 1.2ex}
\hrule}}}

\def \outputWt{L^{-1}}

\usepackage{color}

\def \figname {fig:minplusrobustest:}
\def \eqnname {eq:minplusrobustest:}

\def \secname {sec:minplusrobustest:}

\def \obsvwindowlength{\mathcal{N}}
\def \disturbanceSet{\mathcal{Z}}

\author{
\authorblockN{Srinivas Sridharan}
\authorblockA{Dept. of Mechanical and Aerospace Engineering \\
University of California San Diego \\
Email: {srsridharan@eng.ucsd.edu}}
\thanks{Research supported by AFOSR Grant FA 9550-10-1-0233.  }
\and
\authorblockN{William M. McEneaney}
\authorblockA{Dept. of Mechanical and Aerospace Engineering \\
University of California San Diego \\
Email: {wmceneaney@ucsd.edu}}
\thanks{Research partially supported by 
NSF grant DMS-0808131 and AFOSR.}
}

\title{Deterministic  filtering and dimensionality reduction for optimal attitude estimation on $SO(3)$} 

\begin{document}
\maketitle
\begin{abstract}
In this article we introduce the use of recently developed min/max-plus techniques 
in order to solve the optimal attitude estimation problem in filtering for nonlinear systems on the  special orthogonal ($SO(3)$) group.
This work helps obtain computationally efficient methods for the synthesis of 
deterministic  filters for nonlinear systems -- i.e. optimal filters which estimate the state using a related
optimal control problem. The technique indicated herein is validated using a set of optimal attitude estimation example 
problems on  $SO(3)$.
\end{abstract}
\section{Introduction}\label{\secname introduction}
Optimal filtering is one of the major themes of research interest
in the areas of systems and control. There are two distinct approaches
to filter design which have been applied in the design of filters for 
systems subjected to process/measurement noise.  The most 
well known approach to filter design is the Kalman filter developed
for linear systems in the 1960\rq s \cite{kalman1960new}. 
The Kalman filter
uses the statistics of the noise/measurement processes in order
to compute the coefficients of the equations used to update the
state estimate as new measurements are observed. An alternative approach
is the minimum energy filtering technique developed by Mortensen \cite{mortensen1968}. This interpretation views the optimal filtering problem as a 
optimal control problem where the objective is to minimize the energy
of the noise process required to explain the observations.
The resulting filter obtained via both these approaches for the 
linear system case with white noise processes (and a quadratic
energy function with $L^2$ noise processes) turns out to be the 
Kalman filter \cite{jazwinski1970stochastic,hijab1980minimum,willems2004deterministic}.

In the case of nonlinear systems (or non-normal noise) the congruence in these two approaches
no longer holds. There are various extensions of the Kalman filter to 
nonlinear systems (e.g. unscented Kalman filter \cite{julier2002reduced,wan2000unscented}, particle
filter \cite{doucet2001sequential}, extended Kalman filters \cite{anderson1979optimal}). However
the extensions of the deterministic filter to the nonlinear 
case and non Euclidean spaces have only recently garnered greater attention \cite{aguiar2006minimum,marcus1984algebraic}. The  specific 
problem of practical  interest  motivating this work on filter 
design for nonlinear systems is the attitude estimation 
problem \cite{bonnabel2009invariant} on the group of rotation matrices (the special orthogonal
group $SO(3)$).

Recent work by Coote et al.\cite{coote2009near} and Zamani et al.\cite{zamani2011near} studied the 
deterministic optimal filtering problem for the case of  systems
evolving on $SO(2)$ and $SO(3)$ respectively. The filter obtained
in the latter case is a near optimal solution to the filtering
problem (i.e.  it is shown to lie within a specific performance
gap from the optimal filter). The motivation in this work on using a near optimal approximation 
has to do with the numerical intractability of solving the optimal control problem associated 
with the filter design. For nonlinear systems standard approaches using the dynamic programming 
method suffer from larger computational requirements (termed the curse of dimensionality) while computing numerical solutions to these classes of problems. 

 In the recent past, a new class of approaches 
to manage this computational intractability  have been developed for both filtering \cite{fleming2000max} and control \cite{mceneaney2008curse,mceneaney2008cdf,SSMc2011Allerton} for a variety of application domains.  
In this article we introduce and solve  the Hamilton-Jacobi-Bellman (HJB)  equation 
for the optimal control problem (associated with the filtering 
problem) using recently developed theoretical and numerical techniques for nonlinear filtering.
These approaches  (termed min/max-plus methods) utilize the fact that the dynamic programming operator for the HJB equation is a linear operator on a particular algebra (the semi-convex functions). Hence, in a certain class of problems, it is possible to solve the HJB equation without having to generate a mesh/grid thereby allowing for numerical tractability.   By applying these techniques to the current problem we indicate an approach to achieve dimensionality reduction in 
filter synthesis for nonlinear systems \cite{sridharan2012min,kallapur2012min} -- specifically the attitude estimation problem; Thus we 
extend the max/min-plus  filtering methods to the case of a compact Lie group 
($SO(3)$).  

This article is organized as follows: we first introduce the system description 
and the problem statement for the attitude estimation problem in 
Sec.~\ref{\secname sysAndProbDescription}.  This leads to a consideration of the variation 
in the parameterization of the value function as it is propagated in time using the dynamic programming propagator. It is shown in Sec.~\ref{\secname minplusandstructpres} that this structure is preserved. Hence this permits the use of a specific expansion of the value function in terms of a basis whose form is preserved -- allowing for the generation of reduced dimensionality
methods for filtering.  The theory developed is then applied to a class of example problems
in Sec.~\ref{\secname example} thereby validating the techniques introduced.  In Sec.~\ref{\secname conclusion} we conclude this article with a description of  several interesting avenues for future research into different aspects of the min-plus approach to filtering for the classes of problems described herein.

\section{System description}\label{\secname sysAndProbDescription}
Let the original system dynamics be 
\begin{align}
\dot{R} = R (A + z), \label{\eqnname gensysmodel} \\
Y = R \epsilon, \qquad \epsilon \in SO(3),
\end{align}
where $Y, A$ are known signals and $z$ is the unknown state disturbance signal and $\epsilon$ is the unknown measurement noise. 
Let the \textit{backward} time state transition function be defined as
\begin{align}
R_k = \psi(R_{k+1},z_{k+1},A_{k+1}) := R_{k+1} \Psi(z_{k+1},A_{k+1}).
\end{align}
Thus the operator $\Psi(\cdot)$ is the state transition operator, which in this case would be time ordered exponential map. 
Given measurements $Y$, drift $A$, and an initial state estimate $\hat{R}_0$ the cost function for the filtering problem is given by \cite{zamani2011near}
\begin{align}
V_0(R) \aeq \inf_{z \in \disturbanceSet} \Big\{ \int_0^T{ \frac{1}{2} \trace(z^Tz) ds} + \ldots \noncr
& \int_0^T{\frac{1}{4} \phi_{\outputWt}(R^{-1}Y) }ds   +
 \ldots \noncr
& \frac{1}{4}\phi_{K^{-1}}(R(0;z, A) - \hat{R}_0) \Big\}. \label{\eqnname valfndefn}
\end{align}
where
\begin{align}
\phi_M(R) := \trace\Big[{(R-I)^T M (R-I)}\Big], \quad R \in SO(3)
\end{align}
and $R(0;z, A)$ denotes the state at time $0$ given the choice  $z$ of the 
disturbance signals which lie in the space $\disturbanceSet$ of $L^2$ signals. Thus the cost function  is a measure of the weighed sum of the energy in the disturbance in the system dynamics, the unexplained part of the measurements and the error in the initial estimate (versus the predicted initial estimate as obtained from the choice of the disturbance signal). Thus the terminal cost at a point $R_0$ with an initial state estimate $\rhatini$ is defined as
\begin{align}
V_0(R_0) \aeq \frac{1}{4}\phi_{K^{-1}}(R_0 \rhatini^T) \noncr
\aeq \frac{1}{4}  \trace\Big[{(R_0 \rhatini^T-I)^T {K^{-1}} (R_0 \rhatini^T-I)}\Big].\label{\eqnname termvalfn1}
\end{align}
Using the orthogonality property of $R_0$ and $\rhatini$
\begin{align}
R_0^T R_0 = I  = R_0 R_0^T, \qquad \rhatini^T \rhatini = I = \rhatini \rhatini^T, 
\end{align}
 the following properties of the trace operator
\begin{align}
\trace[A B ] = \trace[BA], \noncr
\trace[P] = \trace[P^T],
\end{align}
and the symmetric nature of $K$, it follows that \eqref{\eqnname termvalfn1} is  of the form
\begin{align}
 \aeq \frac{1}{2} \trace\Big[ (K^{-1} - K^{-1} R_0 \rhatini^T)\Big],
\end{align}
Hence, the value function has an affine structure  of the form
\begin{align}
V_0(R_0)  \aeq c_0 + P_0(R_0),   \noncr \text{where} \quad c_0 \adefeq(1/2) \trace(K^{-1}),\noncr P_0(R)\adefeq -(1/2)\trace(\rhatini^T K^{-1} R). \label{\eqnname affinevalfn}
\end{align}
The latter equation \eqref{\eqnname affinevalfn} is obtained using the circular property of the trace operator. 
Note that terminal cost function \eqref{\eqnname affinevalfn} is affine in $R_0$\footnote{More precisely, the cost function is afine if $R_0$ is considered to be the element of $SO(3)$ after embedding in a vector space. Note also that the penalty function $\phi_M(R_0 \rhatini^T)$ is equal to the alternative penalty function $\tilde{\phi}_M(R_0, \rhatini):= \trace[(R_0 - \rhatini)^T M (R_0 - \rhatini)]$, where $R_0$ and $\rhatini$ are treated as elements of the vector space.}. 

In order to solve for the value function we apply the dynamic programming method\cite{bellman2003dp} from optimal control theory. For the cost function \eqref{\eqnname valfndefn} the  dynamic programming principle takes the form
\begin{align}
V_s(R) \aeq 
 \inf_{z \in \disturbanceSet[s, s+\tau]} \Big\{ \int_s^{s+t}{ \frac{1}{2} \trace(z^Tz) dt} + \ldots \noncr
& \int_s^{s+\tau}{\frac{1}{4} \phi_{\outputWt}(R^{-1}Y) }dt   +
 \ldots \noncr
& V_{s+\tau}(R(s+\tau; z, A) \Big\}, \label{\eqnname dpeforvalfn}
\end{align}
where for $t_1, t_2 \in \mathbb{R}^{+}$ s.t $t_1<t_2$, $\disturbanceSet[t_1,t_2]$ denotes the restriction of the control signal space to the time horizon $[t_1, t_2]$.

 We 
first indicate the effect of one step of propagation of the value function \eqref{\eqnname valfndefn}, using the dynamic programming equation \eqref{\eqnname dpeforvalfn}, on its affine structure.
In  order to solve for the value function numerically we make the following assumption. We use a 
discretized and bounded version of the disturbance set $\disturbanceSet$. For simplicity of notation 
we abuse the notation and reuse it to denote the bounded discrete set.  Note that the boundedness 
assumption is  not overly restrictive for reasons that will be explained in Sec.~\ref{\secname subsetpropagation}.

Assuming a discretization time $\Delta t$ the application of the dynamic programming principle to the value function yields
\begin{align}
&V_1(R) = \inf_{z \in \cZ}\Big\{ \frac{1}{2} \trace(z^T z) \Delta t + \frac{1}{4}\phi_{\outputWt}(R^{-1}Y) \Delta t + \ldots \noncr &\qquad V_0(R(0;z,A,Y) \Big\},\\
 &= \inf_{z \in \cZ}\Big\{ \frac{1}{2} \trace(z^T z) \Delta t + \frac{1}{4}\phi_{\outputWt}(R^{-1}Y) \Delta t + \ldots \noncr & \qquad \frac{1}{4} \phi_{K^{-1}}(\psi(R, z, A) \rhatini^T)   \Big\},\\
&= \inf_{z \in \cZ} \Big\{ \frac{1}{2} \trace(z^T z) \Delta t +\frac{1}{2} \trace\Big[\outputWt -\outputWt R^{-1} Y \Big] \Delta t +  \ldots \noncr & \qquad \frac{1}{2} \trace\Big[K^{-1} - \rhatini^T K^{-1}R\Psi(A_1,z)  \Big]  \Big\}.
\end{align}
If $z$ belongs to a discretization of the control set $\mathfrak{so}(3)$, then the above expression is affine in $R \in SO(3)$ for each element $z$ in the algebra. It takes the form
\begin{align}
V_1(R) \aeq \inf_{z \in \cZ}\Big\{c_z + P_z(R)\Big\} =  \bigoplus_{z \in \cZ} c_z \otimes P_z(R),  \end{align}
where
\begin{align}
 c_z \adefeq  \frac{1}{2} (z^Tz) \Delta t + \frac{1}{2} \trace(\outputWt)\Delta t + \frac{1}{2} \trace(K^{-1}), \noncr
 P_z(R) \adefeq  -\frac{1}{2} \trace({\outputWt}^T Y^T R  \Delta t+  \rhatini^T K^{-1}R\Psi(A_1,z) ),\noncr 
 \aeq   -\frac{1}{2} \trace \Big( \big[{\outputWt}^T Y^T + \Psi(A_1,z)  \rhatini^T K^{-1}\big] R \Big).
\end{align}
Thus the min-plus affine structure of the value function is preserved after one step (from the initial time step). This motivates the following more general proof of the invariance of this affine property.
under propagation by the dynamic programming operator. 
\section{The min-plus expansion and the propagation of the cost function}\label{\secname minplusandstructpres}
We first indicate the structure preservation and then describe how this approach to 
solving for the optimal cost function can be used to obtain the state estimate subsequent
to propagation of the solution for the duration of the filter time horizon.
\subsection{Propagation of the form of the cost function} \label{\secname subsetpropagation}
Let the form of the max-plus basis expansion of the value function at time step $k$ be 
\begin{align}
V_k(R) \adefeq \bigoplus_{\lambda_k \in \Lambda_k} \Big[ c_{\lambda_k}  \otimes P_{\lambda_k}(R)\Big], \label{\eqnname minplusbasisexpvalfn}
\end{align}
where
$P_{\lambda_k}(R) = -\frac{1}{2}\trace(M_k R)$.
Now, from the dynamic programming principle 
\begin{align}
V_{k+1}(R) \aeq \bigoplus_{z \in \cZ} \Big(  \frac{1}{2} \trace(z^T z) \Delta t + \ldots \noncr &\qquad \frac{1}{2} \trace\Big[\outputWt -\outputWt R^{-1} Y \Big] \Delta t + \ldots  \noncr & \qquad  V_{k}(R\Psi(A_{k+1},z))   \Big),\noncr
\aeq \bigoplus_{\lambda_{k+1} \in \Lambda_{k+1}} \Big[c_{\lambda_{k+1}}+ P_{\lambda_{k+1}}(R)\Big].
\end{align}
where, by comparing coefficients, it is seen that
\begin{align}
c_{\lambda_{k+1}} \aeq c_{\lambda_{k}} + \frac{1}{2} \trace(z^Tz) \Delta t + \frac{1}{2} \trace(\outputWt)\Delta t, \noncr
P_{\lambda_{k+1}}(R)\aeq -\frac{1}{2} \trace(M_{k+1} R ), \end{align}
where
\begin{align}M_{k+1}  = {\outputWt}^T Y^T  \Delta t +  \Psi(A_{k+1},z) M_k.
\end{align}
Note that for a bounded continuous value function on a compact set, the  optimal value
of the disturbance signal $z$ is bounded. This is the reason for the aforementioned non-restrictiveness of the boundedness assumption.
The basis function  parameterization sets  above evolve via the form
$\Lambda_{k+1} = \Lambda_k \times \cZ$. The above expressions provide the mechanism via which the basis function  sequences $\lambda_{(\cdot)}$ and their corresponding components $c_{(\cdot)}$ and $P_{(\cdot)}$  are generated at each time step in response to the new measurements being gathered therein.

\subsection{Generating the state estimate}
Having obtained the value function in terms of  a min-plus basis expansion \eqref{\eqnname minplusbasisexpvalfn}  we have the following expression for the optimal state estimate.
Given the observation time horizon $\obsvwindowlength$ (the window over which we use the observations to generate a state estimate), the optimal state estimate
is defined as follows
\begin{align}
V_{\obsvwindowlength}(R) \adefeq \bigoplus_{\lambda \in \Lambda_\obsvwindowlength} \Big[ c_{\lambda}  \otimes P_{\lambda}({\cdot})\Big](R), 
 \end{align}
 where $\Lambda_\obsvwindowlength$ is the collection of affine functions generated  at the time step $\obsvwindowlength$ via the
 recursive propagation as described in the previous section. 
In previous approaches to min-plus fitlering  for general non-linear systems on Euclidean space, quadratic basis functions
led to direct approaches to obtain this estimate (determining the minimum of each of the quadratic functions yields the desired 
state estimate c.f. \cite{kallapur2012min}). However in the case of system dynamics evolving on manifolds, such as the $SO(3)$ group considered
herein, a more sophisticated approach
is required owing to the structure of this group. This proceeds by formulating the problem as a convex optimization problem 
for which computationally efficient techniques exist \cite{ben2001lectures}. The details of this approach 
are as follows.

The minimum value of the optimal cost (value) function is  achieved at the optimal state estimate. Hence the estimated state is obtained as
\begin{align}
\hat{R}^*= \argmin_{R \in SO(3)} \bigoplus_{\lambda \in \Lambda_\obsvwindowlength} \Big[ c_{\lambda}  \otimes P_{\lambda}({R})\Big]. \label{\eqnname combinedoptimprob} 
\end{align}
The above problem \eqref{\eqnname combinedoptimprob} can be solved faster, by breaking it down into a set of independent parallel tasks, each of which solves a  problem of the form
\begin{align}
\hat{R}^*_\lambda := \argmin_{R \in SO(3)}\Big[ c_{\lambda}  \otimes P_{\lambda}({R})\Big],\quad \lambda \in \Lambda_\obsvwindowlength.
\end{align}
\section{Example} \label{\secname example}
In this section we describe an application of the technique developed
in the previous sections to a sepcific set of examples of a tracking 
problem on $SO(3)$. This application will serve to indicate the 
performance of the proposed approach in the case of the bilinear 
model \eqref{\eqnname gensysmodel} introduced in Sec.~\ref{\secname introduction}.

We use the following notation for the  directions used to generate the disturbance set $\cZ$.  They lie along the positive and negative directions of the skew-symmetric basis elements of $\mathfrak{so}(3)$. 
\begin{align}
H_1 &:= \left(\begin{array}{ccc}0 & +1 & 0 \\-1 & 0 & 0 \\0 & 0 & 0\end{array}\right),  \noncr H_3 &:= \left(\begin{array}{ccc}0 & 0 & +1 \\0 & 0 & 0 \\-1 & 0 & 0\end{array}\right), \noncr H_5&:=\left(\begin{array}{ccc}0 & 0 & 0 \\0 & 0 & +1 \\0 & -1 & 0\end{array}\right), \\
H_2 &:= -H_1, \qquad H_4 := -H_3, \qquad H_6 := -H_5.
\end{align}
The system dynamics considered in this section conform to the general form \eqref{\eqnname gensysmodel} 
and we apply min-plus approach for specific cases of the model in order to better understand the performance of this approach.  The first three cases below start with a zero initial estimation error (however starting from the first time step an error arises which is the reason for the nonzero starting value of the error in the figures).
\begin{enumerate}
\item{\bf Zero drift with collinear  measurement and process disturbance:}  In this case the dynamics
is  $\dot{R} = R z$ where $z$ is generated as a normal distribution along the direction $H_1$. In 
addition the measurement noise are along the same direction $H_1$ (ref Fig.\ref{\figname zeroDriftCase}).
\begin{figure}[htp]
\begin{center}
 \includegraphics[width=\hsize]{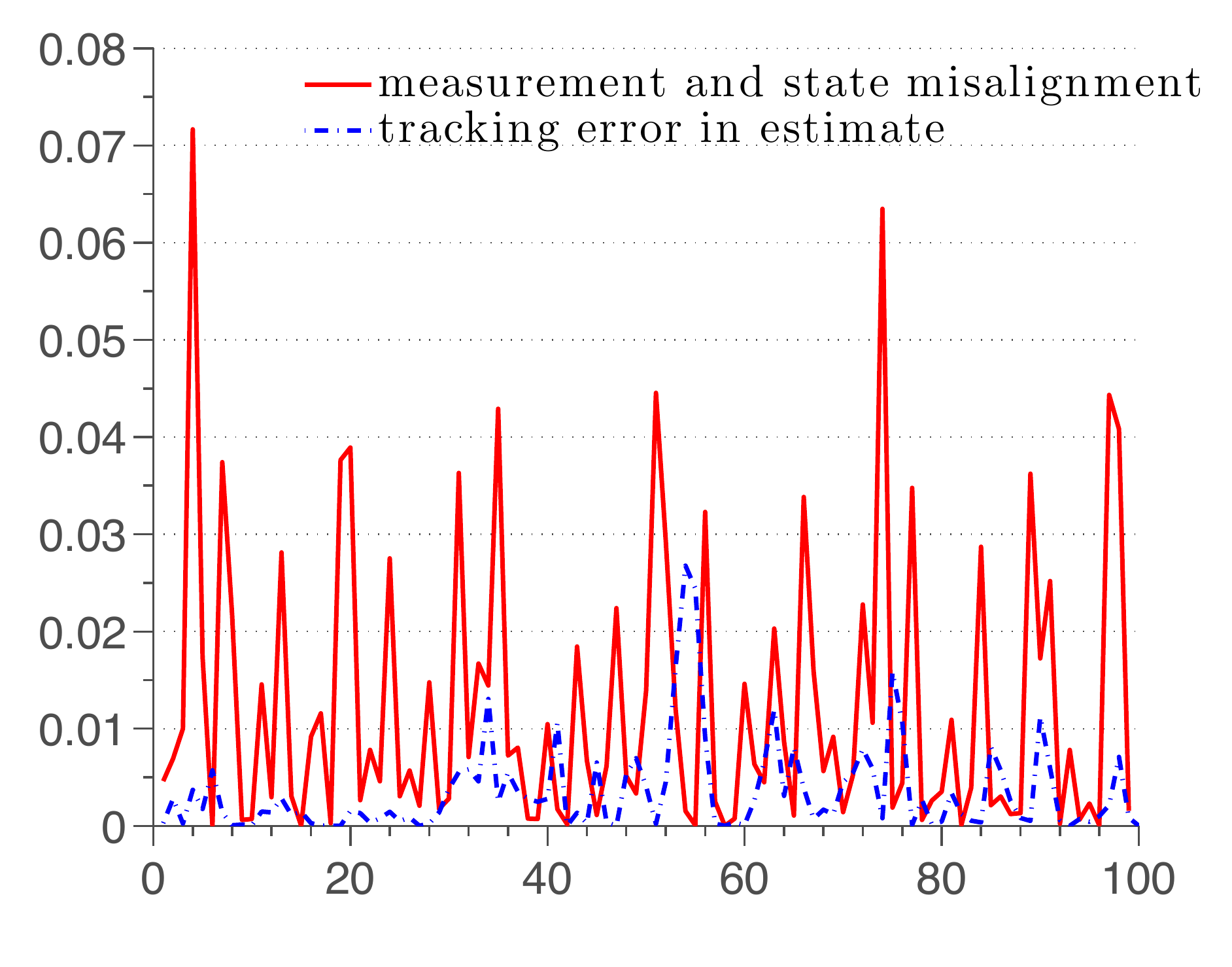}
\caption{The case of zero drift with the dynamic and measurement noise along $H_1, H_2$.} \label{\figname zeroDriftCase}
\end{center}
\end{figure}

\item{\bf Non-zero drift with collinear measurement and process disturbance:} In this case the dynamics 
is $\dot{R} = R (H_1+z)$ and the state disturbance and measurement errors are normally 
distributed along $H_1$ (c.f. Fig.~\ref{\figname nonzerodriftcoaxialMeasAndDisturbance}).

\begin{figure}[htp]
\begin{center}
 \includegraphics[width=\hsize]{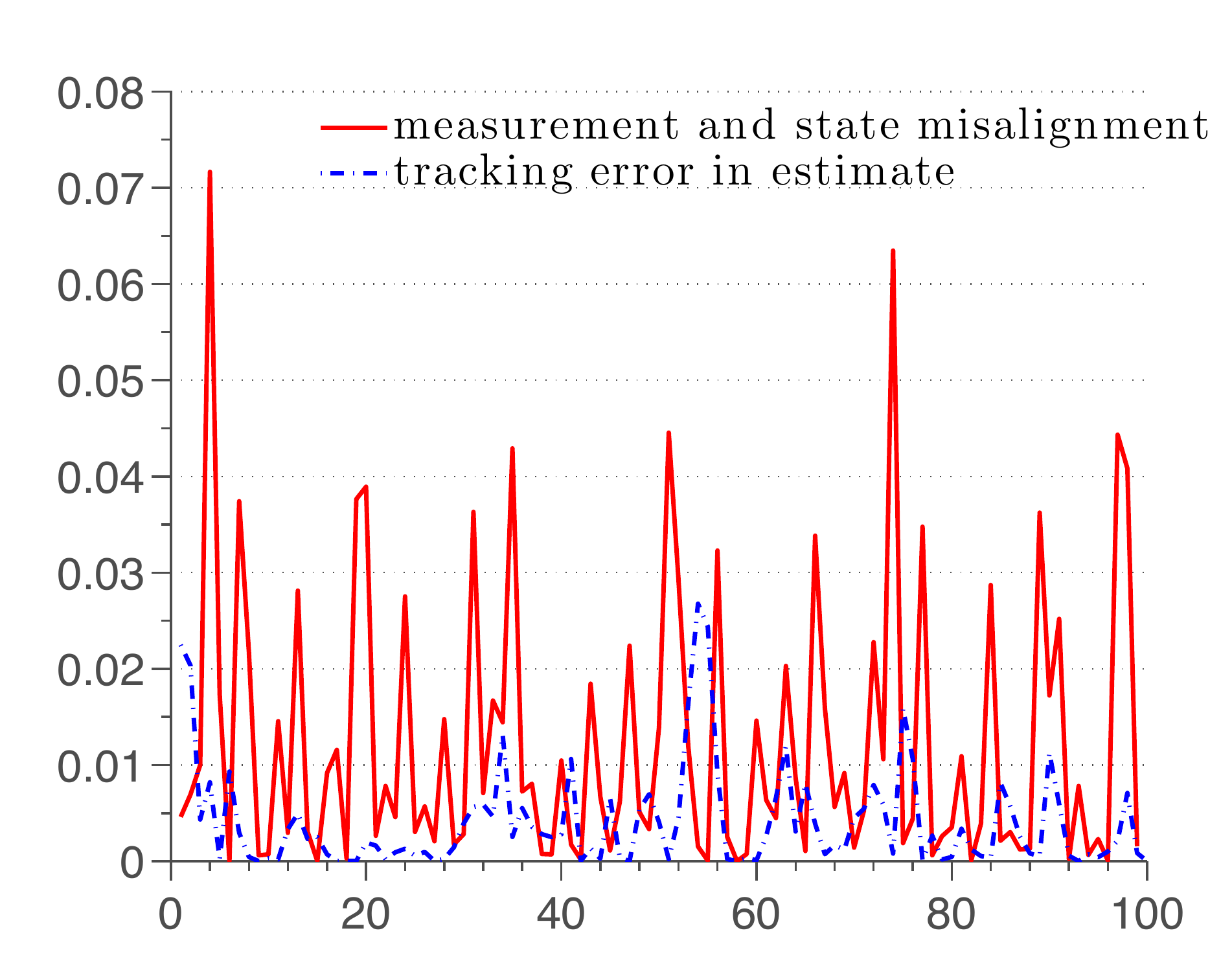}
\caption{The case of  drift along $H_1$ with the dynamic and measurement noise along $H_1, H_2$.} \label{\figname nonzerodriftcoaxialMeasAndDisturbance}
\end{center}
\end{figure}
\item{\bf Non-zero drift with orthogonal measurement and process disturbances:} In this case, the drift direction
is $H_1$ and the measurement noise $\epsilon$ and disturbance $z$ are along the $H_1, H_2$ and $H_3, H_4$ directions respectively (c.f. Fig.~\ref{\figname nonzerodriftNonParallelMeasAndDisturbance}).

\begin{figure}[htp]
\begin{center}
 \includegraphics[width=\hsize]{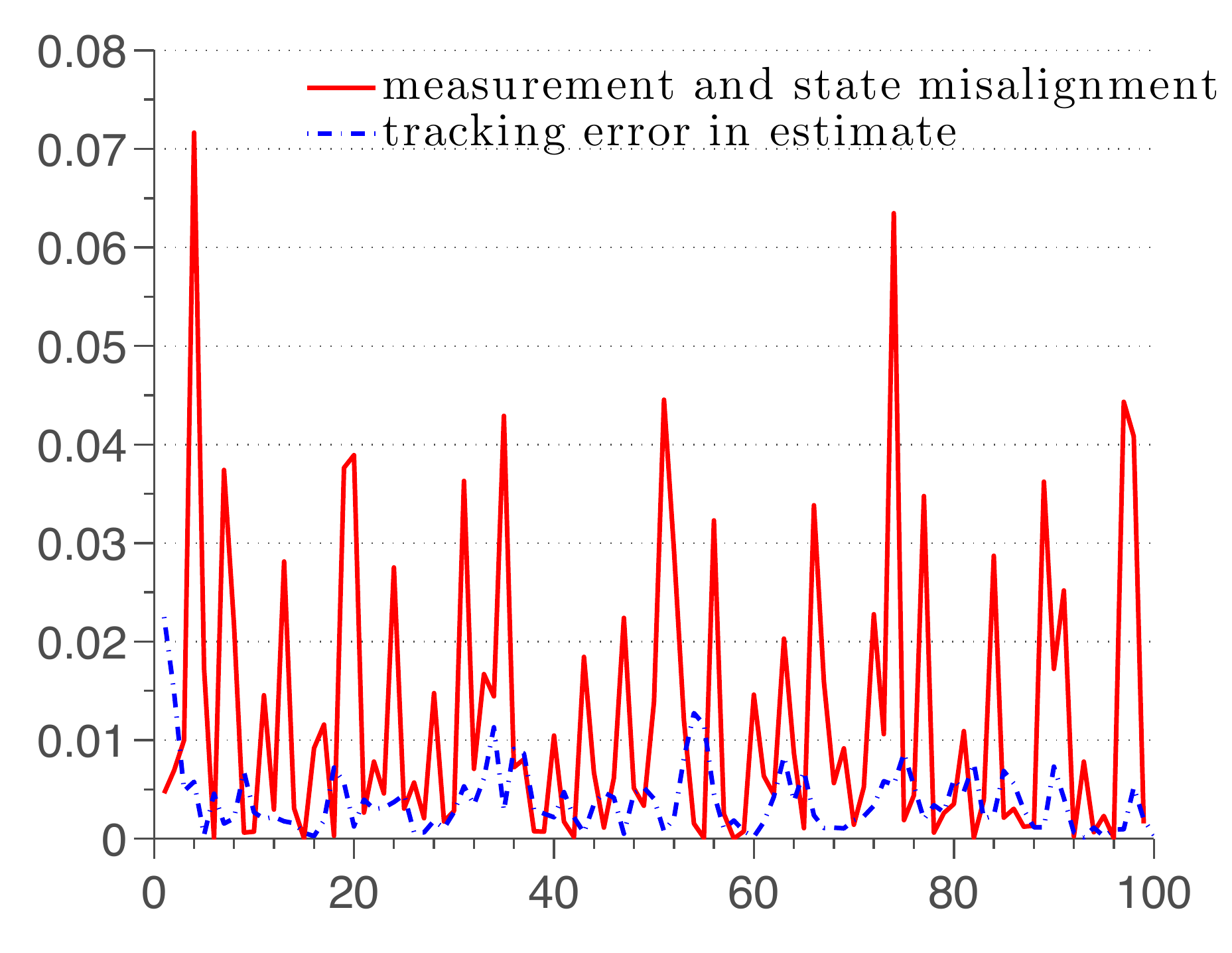}
\caption{The case of  drift along $H_1$ with the dynamic and measurement noise along $H_3,H_4$ and $H_1, H_2$ respectively.} \label{\figname nonzerodriftNonParallelMeasAndDisturbance}
\end{center}
\end{figure}
\item{\bf Non-zero drift with initial estimation error and orthogonal measurement and process disturbances:} In this case, the drift direction
is $H_1$, there is an initial error in the state estimate; the measurement noise $\epsilon$ and disturbance $z$ are along the $H_1, H_2$ and $H_3, H_4$ directions respectively (c.f. Fig.~\ref{\figname nonzerodriftNonParallelMeasAndDisturbance}).

\begin{figure}[htp]
\begin{center}
 \includegraphics[width=\hsize]{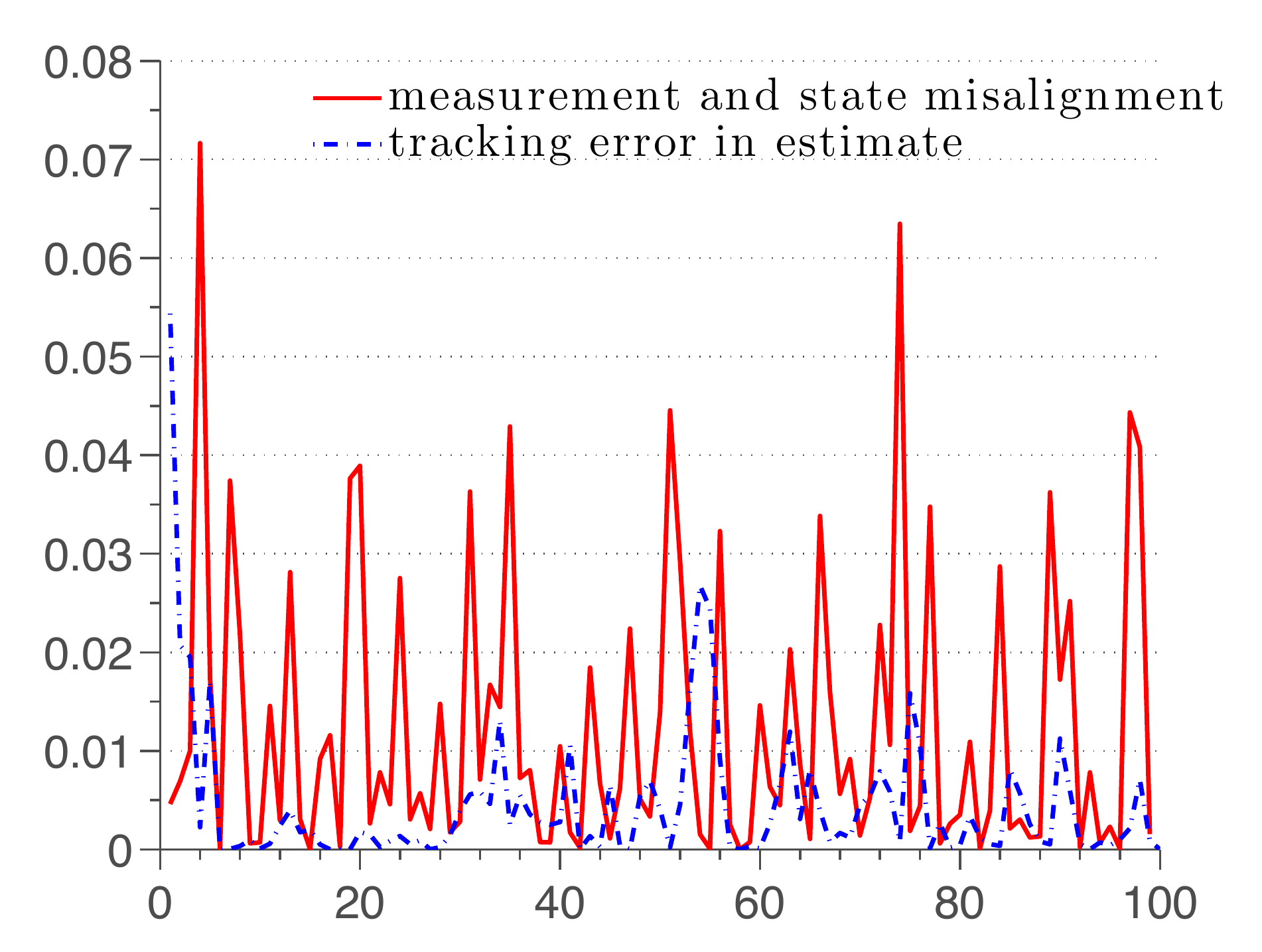}
\caption{The case of  drift along $H_1$ with an initial error in the state estimate; the dynamic and measurement noise  are along $H_3,H_4$ and $H_1, H_2$ respectively.} \label{\figname h1DriftWithh1distNormalInitialErrorh1measNoise}
\end{center}
\end{figure}
\end{enumerate}
In the above models we considered a Gaussian noise characteristic for the process and measurement noise. Below we indicate the performance for the case of noise generated by a 
uniform random process.
\begin{figure}[htp]
\begin{center}
 \includegraphics[width=\hsize]{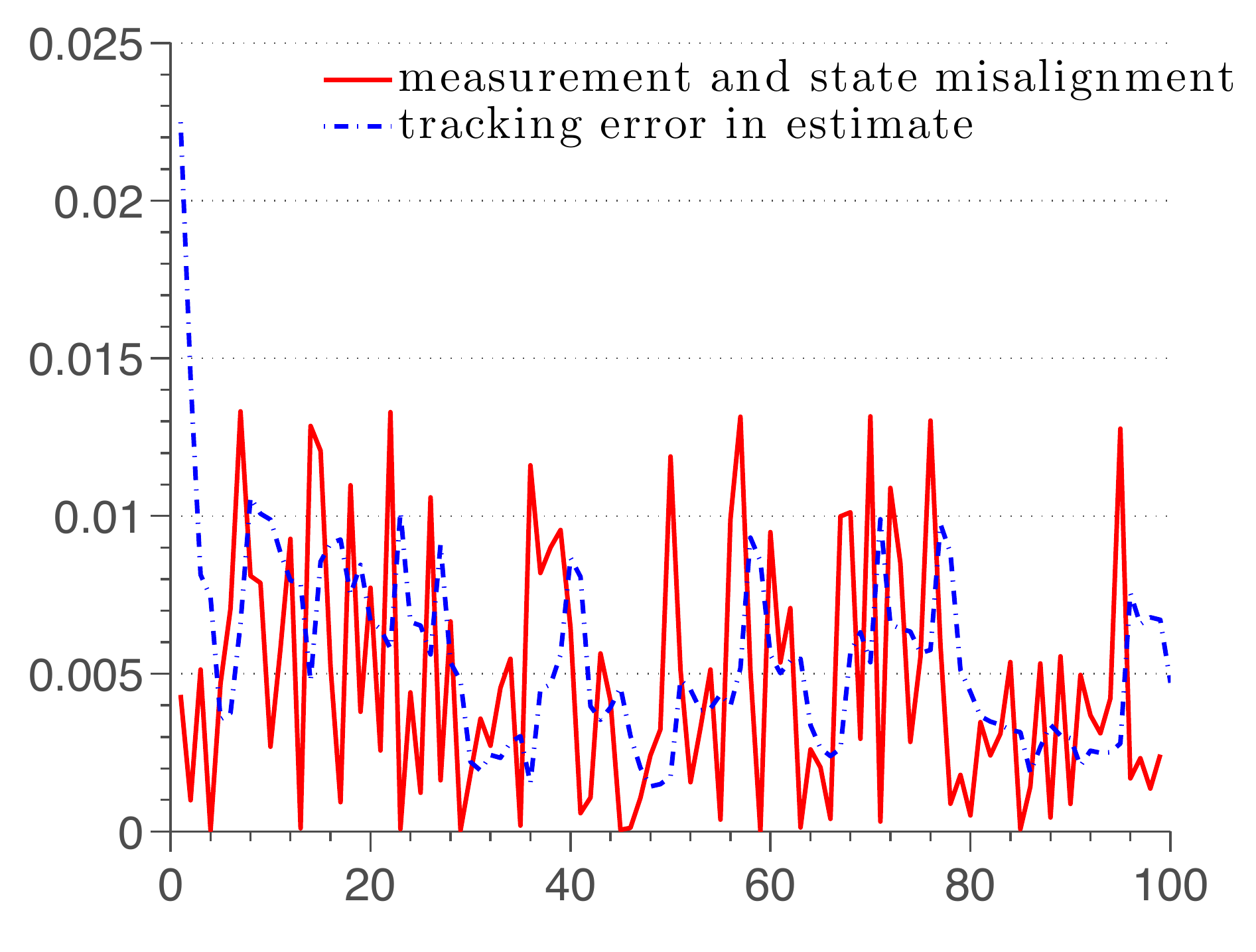}
\caption{The case of  drift along $H_1$ with the dynamic and measurement noise along $H_3,H_4$ and $H_1, H_2$ respectively. In this case the noise samples are from a {\textit{ uniform distribution}}.} \label{\figname nonzerodriftUniformAndNonParallelMeasAndDisturbance}
\end{center}
\end{figure}

Note that in the figures, we analyze the performance of the filter by comparing the measurement noise and the error in the estimates.  The metric which we use   
to indicate the level of noise in the measurements is
\begin{align}
\text{measurement noise} (t) := \trace\Big[I - Y^{T}(t) R(t) \Big],
\end{align}
at a time $t$; the measure of the agreement between the estimate and the true state i.e., the error in the estimate is given by the tracking error in the estimate (TE)
\begin{align}
\text{TE} (t) := \trace\Big[I - {\hat{R}}^{T}(t) R(t) \Big].
\end{align}

In the cases above, we apply the min-plus reduced dimensionality techniques with fixed values of the weighting matrices $(K, L)$ and with 
a fixed length of the time horizon over which the optimization is carried out (width of the sliding window).

Note that the parameters in the studies above were not selected using any optimization criteria. They serve to be indicative 
of the levels of performance achievable with minimal effort and hence demonstrate the potential nominal performance achievable
from these classes of filters with little effort at tuning them. This also serves to guide
some of the directions which such work could proceed along as described in the next section.

\section{Conclusions and Future Directions} \label{\secname conclusion}
In this article we describe the min-plus approaches for deterministic filtering on the SO(3) group. The approach introduced
herein was then applied to example problems for a specific class of dynamics and noise processes. The 
results obtained reveal the need to study various aspects  of these class of techniques. Some of the paths that this work can 
be directed along are: the study of the sensitivity of various filter weight matrices on the performance of the filter; the error 
analysis for a specific value of these parameters; an analysis of 
the length of the filtering window and weight matrices on the performance of the filter; determining the optimal values 
of these filter parameters for a specific noise process having known statistics; analysis of the filter robustness to unknown
noise/disturbance and unmodeled dynamics; the consideration of min-plus based deterministic filter design for alternative system models \cite{bonnabel2009invariant} for attitude tracking problems.
Note that in the current article we assume that the number of terms in the basis expansion generated  for the given filter time 
horizon is capable of being stored in memory. However in the case that this growth in the number of terms if unfeasible, there exist
pruning methods \cite{sridharan2012min} which can be considered in order to manage this growth (e.g. when the time horizon of the filter or the cardinality 
of the disturbance process set is larger than the available memory constraints). Thus this article provides a potential core for several 
 diverse lines of investigation into the area of  computationally efficient approaches to deterministic filtering for nonlinear systems  evolving on manifolds.
\balance
\bibliographystyle{unsrt}

%

\end{document}